\documentclass{amsart}
\usepackage{amsfonts}

\def\RR{\mathbb{R}}
\def\CC{\mathbb{C}}
\def\HH{\mathbb{H}}

\def\OO{\mathbb{O}}

\def\hhh{\mathfrak{h}}

\def\qedbox{\hbox{$\rlap{$\sqcap$}\sqcup$}}
\def\qed{\nobreak\hfill\penalty250 \hbox{}\nobreak\hfill\qedbox}

\def\Im{\mbox{Im}\,}
\def\Id{\mbox{Id}\,}

\newtheorem{Theorem}{Theorem}[section]
\newtheorem{Corrolary}{Corrolary}[Theorem]
\newtheorem{Definition}{Definition}[section]
\newtheorem{Remark}{Remark}[Theorem]
\newtheorem{Lemma}{Lemma}[section]

\begin{document}

\title{Two Classes of Slant Surfaces in Nearly K\" ahler  Six  Sphere}

\author{K. Obrenovi\' c}
\email{kobrenovic@gmail.com}

%    author two information
\author{S. Vukmirovi\' c}
\address{Faculty of mathematics\\ University of Belgrade\\Studentski Trg 16\\11001 Belgrade\\ SERBIA }
%\curraddr{}
\email{vsrdjan@matf.bg.ac.rs}
\thanks{This work was supported by the Ministry of Science, Republic of Serbia, projects 144026 and 144032}



\date{\today}

\dedicatory{}

\begin{abstract}
In this paper we find  examples of  slant surfaces  in the nearly K\" ahler six  sphere.
First, we characterize two-dimensional small and great spheres which are slant. Their description is given in terms of the associative $3$-form in $\Im \OO .$
Later on, we classify the slant surfaces of $S^6$ which are orbits of maximal torus in $G_2.$ We show that these orbits are flat tori which are linearly full in $S^5\subset S^6$ and that their slant angle is between  $\arccos \frac{1}{3}$ and $\frac{\pi}{2}.$  Among them we find one parameter family of minimal orbits.
\end{abstract}

\maketitle

\section{Introduction}

\noindent
It is known that  $S^2$ and $S^6$ are  the only spheres that admit an almost complex structure. The best known Hermitian almost complex structure $J$ on $S^6$ is defined
 using octonionic multiplication. It not integrable,  but satisfies  condition $(\nabla _X J)X = 0,$ for the Levi-Civita connection $\nabla$ and every vector field $X$ on $S^6.$ Therefore, sphere $S^6$ with this structure $J$ is usually referred as nearly K\" ahler six sphere.

\noindent
Submanifolds of  nearly Ka\" hler sphere $S^6$ are subject of intensive research.
A. Gray \cite{Gray} proved that almost complex submanifolds of nearly K\" ahler $S^6$  are necessarily two-dimensional and minimal. In  paper  \cite{Bry} Bryant  showed that any Riemannian surface can be embedded in the six sphere as an almost complex submanifold. Almost complex surfaces were further
investigated in  paper \cite{BVW94} and classified into four types.

\noindent
Totally real submanifolds of $S^6$ can be of dimension two or three. Three-dimensional totally real submanifolds are investigated in  \cite{Ejiri} where N. Ejiri proved that they have to be minimal and orientable.
In  paper \cite{BVW97} authors classify totally real, minimal surfaces with constant curvature. Totaly real and minimal surfaces with Gauss curvature $K\in [0,1]$ (and compact) or those with $K$ constant must have either $K=0$ or $K = 1$ (see \cite{DOVV}).

\noindent
Slant submanifold $N$  of an almost Hermitian manifold $(M,g,J)$ is a generalization of totaly real and almost complex submanifold. The latter have  slant angle $\frac{\pi}{2}$ and $0,$ respectively. Slant submanifolds with slant angle $\theta \in (0, \frac{\pi}{2})$ are called proper slant submanifolds. Note that the notion of surface with constant K\" ahler angle coincides with the notion of  slant surface  and in this paper we use the latter one.
General theory regarding slant submanifolds and some classification theorems of slant surfaces in $\CC ^2$ can be found in \cite{Chen1}.
Because of dimensional reasons (see \cite{Chen1}) a proper slant submanifold of  six sphere  is two-dimensional. Very few examples of slant surfaces of almost K\" ahler six sphere are known and all of them are minimal. In paper \cite{BVW94} it is shown that rotation of an almost complex curve of type (III) results in minimal slant surface that is linearly full in some $S^5\subset S^6$. According to \cite{XingMin}, a minimal slant surface of $S^6$ that has non-negative Gauss curvature $K\geq 0$  must have either  $K\equiv 0$ or $K\equiv 1.$ Classification of such surfaces is given in paper  \cite{Xing}. In the present paper we find two classes of slant surfaces of $S^6$ with $K\equiv 0$ and $K\equiv 1.$ Some of them are minimal and therefore known, but most of them are not minimal.

\noindent
In the Section \ref{sec:prelims} we recall some basic facts about octonions and define the notion of a slant submanifold (Definition \ref{de:slant}). Lemma \ref{the:nezavodvekt} is simple, but it is not  found in \cite{Chen1}.
In  Section  \ref{se:spheres} we consider  slant two-dimensional spheres which are intersection of an affine $3$-plane and the six sphere. We give their
characterization in terms of associative $3$-form in Theorem \ref{th:preseci}.
Finally, in the Section \ref{se:orbits} we investigate two-dimensional orbits of a Cartan subgroup of group $G_2$ on the sphere $S^6.$  We found that these orbits are flat, two-dimensional tori which are always slant, but it is interesting that their slant angle is in the range $[\arccos \frac{1}{3}, \frac{\pi}{2}]$ (see Theorem  \ref{th:orbits}). In the Theorem \ref{th:minimal} we find one-parameter family of minimal orbits. Note that similar method was used in paper \cite{HaMa} for obtaining  three-dimensional orbits which are CR submanifolds of sphere $S^6.$

\section{Preliminaries}
\label{sec:prelims}

Let  $\HH$ be a field of quaternions. The Cayley algebra, or algebra of octonions, is vector space  $\OO = \HH \oplus \HH \cong \RR ^8$ with multiplication defined in terms of quaternionic multiplication:
$$
(q,r)\cdot(s,t):=(qs - \overline{t}r,tq+r\overline{s}),\enskip
q,s,r,t\in \HH.
$$
In the sequel we omit the multiplication sign.
Conjugation of Cayley numbers is defined by
$$
\overline{(q,r)}:=(\overline{q},-r), \enskip q,r\in \HH ,
$$
and inner product by:
\begin{equation}
\label{eq:inner}
\langle x,y \rangle := \frac{1}{2}(x  \overline{y} + y
\overline{x}), \enskip x,y\in \OO .
\end{equation}
If we denote by $1,i,j,k$ the standard orthonormal basis of  $\HH,$ then
$e_{0}=(1,0), e_{1}=(i,0), e_{2}=(j,0),
e_{3}=(k,0), e_{4}=(0,1), e_{5}=(0,i),
e_{6}=(0,j), e_{7}=(0,k)$  is orthonormal basis of  $\OO .$ It is easy to check that the multiplication in that basis is given by the following table.
$$
\begin{tabular}{|c|c|c|c|c|c|c|c|}
  \hline
   $ $ & $e_1$ & $e_2$ & $e_3$ & $e_4$ & $e_5$ & $e_6$ & $e_7$ \\
  \hline
  $e_1$ & $-e_0$ & $e_3$ & $-e_2$ & $e_5$ & $-e_4$ & $-e_7$ & $e_6$ \\
  $e_2$ & $-e_3$ & $-e_0$ & $e_1$ & $e_6$ & $e_7$ & $-e_4$ & $-e_5$ \\
  $e_3$ & $e_2$ & $-e_1$ & $-e_0$ & $e_7$ & $-e_6$ & $e_5$ & $-e_4$ \\
  $e_4$ & $-e_5$ & $-e_6$ & $-e_7$ & $-e_0$ & $e_1$ & $e_2$ & $e_3$ \\
  $e_5$ & $e_4$ & $-e_7$ & $e_6$ & $-e_1$ & $-e_0$ & $-e_3$ & $e_2$ \\
  $e_6$ & $e_7$ & $e_4$ & $-e_5$ & $-e_2$ & $e_3$ & $-e_0$ & $-e_1$ \\
  $e_7$ & $-e_6$ & $e_5$ & $e_4$ & $-e_3$ & $-e_2$ & $e_1$ & $-e_0$ \\
  \hline
\end{tabular}
$$
The Cayley numbers are not associative, so we define {\it associator}
$$
[x, y, z]:= (x  y)  z - x  (y  z), \enskip x, y, z \in  \OO .
$$

\noindent
Denote by
$$
\Im \OO := \{ x\in \OO \ |\ x + \overline{x} = 0 \},
$$
the subspace of imaginary Cayley numbers. Then we have the orthogonal decomposition
\begin{equation}
\label{eq:Re+Im}
\OO = \RR \oplus \Im \OO = \RR \oplus \RR ^7.
\end{equation}

\noindent
On the subspace of imaginary Cayley numbers $\Im \OO $ we define the vector product
$$
x\times y := \frac{1}{2}(x  y - y  x),
$$
that shares many properties with the vector product in $\RR ^3.$

\noindent
We state some well known properties of Cayley numbers  without a proof.
\begin{Lemma}
\label{le:cayleyprop}
\begin{itemize}
\item[1)]  If  $x,y\in \Im \OO $ then
$$
x y=-\langle x,y\rangle + x \times y.
$$

\item[2)]  For all  $x,y,z \in \OO$ we have
$$
\overline{x}   (x   y) = (\overline{x}   x)  y,
$$
$$
\langle x  y, x   z\rangle = \langle x,x \rangle \langle
y,z \rangle = \langle y x, z x \rangle.
$$

\item[3)] If $x,y,z \in \OO$ are mutually orthogonal unit vectors then
$$
x  (y  z) = y  (z  x) = z   (x  y).
$$
\end{itemize}
\end{Lemma}

\noindent
Exceptional group $G_2$ is usually defined as a group of automorphisms of Cayley numbers. Since it preserves the multiplication it also preserves the inner product (\ref{eq:inner}) and the decomposition (\ref{eq:Re+Im}) and therefore it is subgroup of group $O(7).$ Actually, it is subgroup of group $SO(7).$

\noindent
For any point $p\in S^6\subset \Im \OO$ and a tangent vector $X\in T_pS^6$ we define automorphism $J_p:T_pS^6\to T_pS^6$ by
$$
J_p(X):= p\cdot X= p \times X.
$$
One can easily show that six-dimensional sphere $(S^6, \langle \rangle, J)$ is almost Hermitian  manifold, i.e. $J_p$ satisfies
$$
J_p^2  = -\Id ,\qquad \langle J_pX, J_pY\rangle = \langle X, Y\rangle,
$$
for all $X, Y\in T_pS^6 .$
\nonumber
The unit six-dimensional sphere  $S^{6}\subset \Im \OO \cong  \RR^{7}$ posses almost complex structure  $J$ defined by
$$
J_{p}(X)=p  X = p\times X,\quad p\in S^{6},  X\in T_{p}S^{6}.
$$

\noindent
Obviously, the group $G_2$ preserves the structure  $J.$

\begin{Definition}
\label{de:slant}
Let  $(M,g,J)$ be an almost Hermitian manifold and   $N\subset M$ be a submanifold of  $M.$
For each $p\in N$ and $X\in T_{p} N$ we define {\it Wirtinger angle}
$$
\theta_{p}(X):=\angle(JX,T_{p}N).
$$
We say that  $N$ is {\em slant submanifold} if its Wirtinger angle $\theta$ is constant, i.e. it doesn't depend on the point  $p\in N$ and tangent vector $X\in T_{p} N$.
\end{Definition}

\noindent
Slant submanifold with angle  $\theta \equiv 0$ is usually called {\em almost complex submanifold}, and slant submanifold with angle $\theta \equiv  \frac{\pi}{2}$ is called {\em totally real submanifold.} Slant submanifold that is neither almost complex nor totally real is called {\it proper slant submanifold.}

\noindent
It is known (see \cite{Chen1}) that proper slant submanifold $N\subset M$ has to be of even dimension and its dimension is less that half of  the dimension of $M.$
The next lemma shows that in the case $\dim N = 2$, the Wirtinger angle of $N$ is always independent on vector $X\in T_pM.$

\begin{Lemma}
\label{the:nezavodvekt}
Let  $(M, g, J)$ be an almost-Hermitian manifold and  $N\subset M$ surface. The Wirtinger angle  $\theta_{p}(X)$ doesn't depend on the vector  $X\in
T_{p}N$ and
$$
\theta_{p}(Z) = |g(X, JY)|
$$
for all $Z\in T_pN,$ where $(X,Y)$ is any orthonormal basis of $T_pN.$
\end{Lemma}

\noindent Proof:
For each $Z \in T_pN$ there is an orthogonal decomposition
$$
JZ = PZ + FZ,
$$
to the tangent component $PZ \in T_pN$ and normal component $FZ  \in  T_p^\perp N.$

\noindent
Therefore,
$$
\cos{\theta_{p}(Z)} = \frac{g(JZ, PZ)}{\|JZ\| \|PZ\|}=
\frac{\|PZ\|^2}{\|Z\| \|PZ\|} = \frac{\|PZ\|}{\|Z\|} .
$$

\noindent
Since $J$ is an isometry we have
$$
g(JX, Y) = -g(X, JY),
$$
for all $X,Y \in T_pM$ and particularly  for $X=Y$
$$
g(JX, X) = 0.
$$
Let $(X, Y)$ be an orthonormal basis of $T_pN.$  For any  $Z= aX +
bY \in T_pN$
\begin{eqnarray*}
 \|PZ\|^2  & = & g(PZ, X) ^2 + g(PZ, Y)^2 = g(JZ, X) ^2 + g(JZ, Y)^2 = \\
& = &  g(aJX + bJY, X) ^2 + g(aJX +
bJY, Y)^2 = \\
& = &  (a^2 + b^2)g(X, JY) ^2,
\end{eqnarray*}
so we have
$$
\cos{\theta_{p}(Z)} = \frac {\sqrt{a^2 + b^2}\, |g(X, JY)|}{\sqrt{a^2+ b^2}} = |g(X, JY)|.
$$
 \qed

\section{Slant two-dimensional spheres in  $S^6$}
\label{se:spheres}

\begin{Definition}
If $\pi\in G_{\RR}(3,\Im \OO)$ is imaginary part or quaternionic subalgebra of Cayley numbers $\OO$ we call
$\pi$ associative $3$-plane. Denote by  $ASSOC\subset  G_{\RR}(3,\Im \OO)$ set of all associative planes.
\end{Definition}
Since the quaternions are associative, the associator of any three vectors of an associative plane is equal zero. Vice versa, if associator of some three vectors from $\Im \OO$ vanishes, then these three elements span an associative plane.

\noindent
On the vector space $\Im \OO = \RR^{7}$ we define $3$-form $\phi$
by the formula
$$
\phi (x,y,z) := \langle x,yz\rangle.
$$
In (\cite{Har}) this form is called  {\it associative  $3$-form} and it is shown that the form $\phi$ is calibration with contact set  $ASSOC.$

\noindent
We use the following notation for the associative $3$-form  $\phi$ and the associator of a $3$-dimensional plane  $\pi \in \Im \OO$
$$
\phi (\pi) :=  |\phi (f_1,f_2,f_3)|, \quad \mbox{and}\quad [\pi] := [f_1,f_2,f_3],
$$
where   $f_1,f_2,f_3$
is orthonormal basis of  $\pi .$
One can show that these definitions do not depend on the choice of orthonormal basis $f_1, f_2, f_3$ of the plane $\pi.$ Furthermore, both the form $\phi$ and the associator are $G_2$ invariant.

\noindent
From the above definition it follows that $\phi (\pi )\in[0,1]$ for any $3$-dimensional plane $\phi$ and that associative planes are characterized by   the condition $\phi (\pi) =  1.$

\noindent
The associator and the associative $3$-form $\phi$ are related by the formula
\begin{equation}
\label{eq:sincos}
\phi^2(\pi) +  \frac{1}{4}\|[\pi]\|^2 = 1 ,
\end{equation}
which we prove in the  Lemma \ref{le:Gram}.

\noindent
In the remainder of this section we characterize two-dimensional spheres that are intersection of $3$-dimensional affine plane and the sphere  $S^6.$
For the proof of Theorem \ref{th:preseci} we need the following lemmas.

\begin{Lemma}
\label{le:Gram}
Let  $f_1, f_2, f_3$ be an orthonormal basis of the plane  $\pi.$ The Gram matrix of the set of vectors
$$
f = (f_1, f_2, f_3, f_2f_3, f_3f_1, f_1f_2, [f_1, f_2,f_3])
$$
is the matrix
$$
\left(
  \begin{array}{ccccccc}
    1 & 0 & 0 & \phi & 0 & 0 & 0 \\
    0 & 1 & 0 & 0 & \phi & 0 & 0 \\
    0 & 0 & 1 & 0 & 0 & \phi & 0 \\
    \phi & 0 & 0 & 1 & 0 & 0 & 0 \\
    0 & \phi & 0 & 0 & 1 & 0 & 0 \\
    0 & 0 & \phi & 0 & 0 & 1 & 0 \\
    0 & 0 & 0 & 0 & 0 & 0 & 4(1-\phi  ^2) \\
  \end{array}
\right),
$$
where we abbreviate $\phi = \phi (f_1, f_2, f_3).$
Particularly, the set $f$ spans  $Im \OO$ if and only if  $\phi (\pi) \neq 1$, i.e. the plane $\pi$ is not associative.
\end{Lemma}
Proof: Most of the scalar product are simple to calculate using properties of Cayley numbers from Lemma \ref{le:cayleyprop}. For example:
$$
\langle f_2f_3, f_3f_1 \rangle = - \langle f_2f_3, f_1f_3\rangle = \langle f_2, f_1\rangle |f_3|^2  = 0.
$$
Now we prove the most complicate inner product $\langle [f_1, f_2,f_3], [f_1, f_2,f_3]\rangle.$
First, note that all $f_2f_3, f_3f_1, f_1f_2$ are imaginary. Using simple transformations one can show
\begin{equation}
f_3(f_1f_2) = -2\phi (f_1, f_2, f_3) - (f_1f_2)f_3.
\end{equation}
Using the previous formula we get
\begin{eqnarray*}
\langle [f_1, f_2,f_3], [f_1, f_2,f_3]\rangle & = &
4(\langle (f_1f_2)f_3 , (f_1f_2)f_3\rangle + 2  \langle (f_1f_2)f_3 , \phi\rangle + \langle \phi, \phi \rangle ) = \\
&=& 4(1-2 \phi^2 + \phi ^2) = 4(1- \phi^2).
\end{eqnarray*}
\qed

\begin{Theorem}
\label{th:preseci}
Let $S^2_r = \pi{'} \cap S^6, \enskip r \in (0,1]$ be a two-dimensional sphere of radius $r.$ Denote by  $\pi \in G_{\RR}(3,\Im \OO)$ $3$-dimensional plane parallel to the affine plane $\pi '$ containing $S^2_r.$
\begin{itemize}
\item[a)]  If $S^2_1$ is great sphere, that is $\pi = \pi^{'},$ then it is slant with the slant angle $\theta = \arccos \phi (\pi)$. Therefore, $S^2_1$ is proper slant sphere if  $\pi$ is not associative.
\item[b)] If the plane $\pi$ is associative  then the small sphere  $S^2_r$ is slant with the slant angle  $\theta = \arccos r$.
\item[c)] If the plane  $\pi$ is not associative  then the small sphere $S^2_r$ is slant with the slant angle  $\theta = \arccos (r\, \phi(\pi))$ if and only if its center is point  $C = \pm \frac{\sqrt{1-r^2}}{|[\pi]\|}[\pi]$.
\end{itemize}
\end{Theorem}

\noindent Proof:  Let $X, Y$ be an orthonormal basis of the tangent space  $T_pS^2_r,$ $p\in S^2_r.$ According to the Lemma  \ref{the:nezavodvekt} we have
$$
\cos\theta _p = |\langle X,JY\rangle|=|\langle X,pY\rangle| =
|\phi(p,X,Y)|.
$$
Particularly, if $p \in S^2_1$ then $p,X,Y$ is an orthonormal basis of the plane  $\pi = \pi '$ and we have
$$
\cos\theta _p = |\phi(p,X,Y)| = \phi(\pi)
$$
which proves the statement a).

\noindent
Now we suppose that $\pi \neq \pi '.$ Let $f_1,f_2,f_3$ be an orthonormal basis of  $\pi .$ Since  $\pi '$ is parallel to  $\pi$ we can write  $\pi '=\pi + \sqrt{1-r^2} \, \xi$ for a unit vector  $\xi \in \pi^{\bot}.$ Therefore $p \in S^2_r$ is of the form
$$
p=rp_0 + \sqrt{1-r^2}\, \xi
$$
for some point   $p_0 = p_1f_1 + p_2f_2 + p_3f_3 \in S^2_1.$ As before, let $X, Y$ be an orthonormal  basis of the tangent space $T_pS^2_r$.
Then vectors $p_0, X, Y$ form an orthonormal basis of  $\pi .$
\begin{eqnarray*}
\cos\theta_p &=& |\phi(p,X,Y)|= |\phi(r\, p_0 + \sqrt{1-r^2}\, \xi, X,Y)|=\\
&=& |r\, \phi(p_0, X, Y)+\sqrt{1-r^2}\, \phi(\xi,X ,Y)|.
\end{eqnarray*}
Since  $|\phi(p_0, X, Y)|= \phi(\pi),$ it remains to calculate
$\phi(\xi, X, Y)$. Note that if  the plane  $\pi$ is associative then  $XY \in \pi$ and we have
$$
\phi(\xi, X, Y)=\langle \xi, XY \rangle = 0.
$$
Therefore, in the case of the associative plane $\pi$ we have
$$
\cos\theta_p = r\phi(\pi) = r,
$$
what proves the statement b).

\noindent
Let  $\pi\neq \pi '$ be a non-associative plane.  Since  $X, Y \in \pi$
we can write them in the form
$$
X= x_1f_1+x_2f_2+x_3f_3,
$$
$$
Y=y_1f_1+y_2f_2+y_3f_3.
$$
The vectors   $X$ and $Y$ are orthogonal to the point  $p_0\in \pi$ its coordinates are
$$
p_0=(p_1, p_2, p_3) = (x_2y_3-x_3y_2, x_3y_1-x_1y_3, x_1y_2-x_2y_1).
$$
Now, easy calculation yields
\begin{eqnarray*}
\phi(\xi, X, Y)& = &\langle \xi,XY\rangle=\langle \xi,
(x_1f_1+x_2f_2+x_3f_3)(y_1f_1+y_2f_2+y_3f_3)\rangle =\\ & = &
\langle \xi, f_2f_3p_1-f_1f_3p_2+f_1f_2p_3 \rangle =\\ & = &
\phi(\xi, f_2, f_3)p_1-\phi(\xi, f_1, f_3)p_2+\phi(\xi, f_1,
f_2)p_3.
\end{eqnarray*}

\noindent
The sphere $S^2_r$ is slant if this expression doesn't depend on point $p\in S^2_r,$ that is, if and only if
$$
\phi(\xi, f_2, f_3)=0, \quad \phi(\xi, f_1, f_3)=0, \quad
\phi(\xi, f_1, f_2)=0.
$$
This condition means that the vector  $\xi$ is orthogonal to vectors  $f_2f_3, f_3f_1, f_1 f_2$, and since $\xi \in \pi ^\perp$, according to the Lemma  \ref{le:Gram}, the only possibility for the unit vector $\xi$ is
$$
\xi = \pm \frac{[\pi]}{|[\pi]|}
$$
and therefore the statement c) holds.\qed

\begin{Lemma}
\label{le:equiv}
Two $3$-dimensional planes  $\pi _1, \pi _2 \in G_{\RR}(3,\Im \OO)$ are $G_2$ equivalent if and only if  $\phi (\pi_1) = \phi (\pi _2)$.
\end{Lemma}
Proof: If the planes $\pi _1$ and $\pi _2$ are equivalent by a $G_2$ transformation  then we have  $\phi (\pi_1) = \phi (\pi _2)$ because the form $\phi$ is $G_2$ invariant. Lets prove the converse.

\noindent
If $\phi (\pi_1) = 1 = \phi (\pi _2)$, i.e. the planes are associative, they are $G_2$ equivalent. Namely, if the planes $\pi _1$ and $\pi _2$ are spanned by orthonormal bases $f_1, f_2, f_1f_2$ and $g_1, g_2, g_1g_2$ respectively, than any $G_2$ transformation that maps  $f_1, f_2$  to $g_1, g_2$ also maps $\pi_1$ onto $\pi_2$.

\noindent
 Suppose that $\phi (\pi) = \phi  \neq 1$ and $f_1, f_2, f_3$ is orthonormal basis of $\pi.$ Denote by $F_1 = f_1, F_2 = f_2, F_3 = f_1f_2$. From Lemma \ref{le:Gram} it follows that  $F_4 = \frac{1}{2\sqrt{1-\phi^2}}[f_1, f_2, f_3]$ is unit vector orthogonal to $F_1, F_2$ and $F_3.$ Following the Cayley-Dixon process the set of vectors  $F_1, F_2, F_3, F_4$, $F_5 = F_1F_4$, $F_6 = F_2F_4$, $F_7 = F_3F_4$ is $G_2$ basis of $\Im \OO$, i.e. satisfies the same multiplication properties as standard basis $e_1, \dots , e_7$ from Section \ref{sec:prelims}.
One can easily check the following relations
\begin{eqnarray}
F_5 &=& \frac{1}{\sqrt{1-\phi ^2}}(\phi f_1 + f_2f_3),\nonumber\\
F_6 &=&  \frac{1}{\sqrt{1-\phi ^2}}(\phi f_2 + f_3f_1),\nonumber\\
F_7 &=& \frac{1}{\sqrt{1-\phi ^2}}(-f_3  + \phi f_1f_2) \label{eq:f3}.
\end{eqnarray}
There is a   $G_2$ transformation that maps vectors  $F_1, F_2, F_4$ to vectors $e_1, e_2, e_4$, respectively. According to the relation (\ref{eq:f3}) the image of vector  $f_3$ is  $\phi \, e_3 -\sqrt{1-\phi^2}\, e_7$.
Therefore, in the case  $\phi (\pi) = \phi  \neq 1$, plane  $\pi$ is   $G_2$ equivalent to the plane
$$
\pi_0^\phi = \RR \langle e_1, e_2, \phi \, e_3 -\sqrt{1-\phi^2}\, e_7\rangle .
$$
We conclude that any two planes  $\pi _1, \pi _2$ with  $\phi (\pi_1) = \phi (\pi _2)$ are $G_2$ equivalent, as claimed.\qed

\begin{Corrolary}
 Let $\pi ' \cap S^6$   and  $\tau ' \cap S^6$ be two-dimensional slant spheres from the Theorem \ref{th:preseci}. They are  $G_2$ equivalent if and only if they are of the same radius and  $\phi (\pi) = \phi (\tau)$.
\end{Corrolary}

\begin{Remark}
Up to $G_2$ equivalence there is only one almost complex two-dimensional sphere. It is great sphere belonging to  associative plane $\pi,$ i.e $\phi (\pi)  = 1.$
\end{Remark}
\begin{Remark}
Totally real  two-dimensional spheres are  $S^2_r = \pi ' \cap S^6$ where $\pi'$ is affine $3$-plane parallel to the plane $\pi$ with $\phi (\pi) = 0$. Up to a $G_2$ equivalence there is unique  such sphere for each radius $r\in (0,1].$ The one with radius $r = 1$ is minimal. It is found in classification of totally real minimal surfaces of constant curvature (see \cite{BVW97}, Theorem 6.5 (a)).
\end{Remark}
\begin{Remark}
Great spheres $S^2_1$ (Theorem \ref{th:preseci} a)) are exactly those from the classification of  minimal slant surfaces with $K\equiv 1$ in $S^6$ (see \cite{Xing}, Example 3.1).
\end{Remark}

\section{Two-dimensional slant orbits in $S^6$  }
\label{se:orbits}

In this section we consider orbits of two-dimensional subgroup  $H\subset G_2\subset SO(7)$ under the natural action on $S^6.$ Since such action preserves both metric on $S^6$ and its almost complex structure all points on a fixed orbit  have the same Wirtinger angle. Therefore, all such two-dimensional orbits are slant surfaces of  $S^6$.

\noindent
One can show that  two-dimensional subgroup of $G_2$ is its Cartan subgroup, i.e. maximal tori of $G_2.$
Since any two Cartan subgroups of $G_2$ are conjugate by some element of $G_2$, they have the same set of orbits. Therefore it is not a loss of generality if we pick any particular Cartan subgroup $H\subset G_2$ to work with.

\noindent
Denote by $E_{[i,j]}=\frac{E_{ij} -
E_{ji}}{2}$,\  $i,j=1,2,...,7,\ i<j,$ the standard basis of Lie algebra  $\mathfrak{so}(7)$ of $SO(7)$. Then a basis of Lie algebra  $\mathfrak{g}_2$ of Lie group $G_2$ is

$$
\begin{array}{cccc}
P_0 = E_{[3,2]}+E_{[6,7]},\qquad & \qquad Q_0 = E_{[4,5]} + E_{[6,7]}, \\
P_1 = E_{[1,3]}+E_{[5,7]}, \qquad& \qquad Q_1 = E_{[6,4]} +
E_{[5,7]},\\
P_2 = E_{[2,1]}+E_{[7,4]},\qquad & \qquad Q_2 = E_{[6,5]} + E_{[7,4]}, \\
P_3 = E_{[1,4]}+E_{[7,2]},\qquad & \qquad Q_3 = E_{[3,6]} + E_{[7,2]}, \\
P_4 = E_{[5,1]}+E_{[3,7]},\qquad & \qquad Q_4 = E_{[2,6]} + E_{[3,7]}, \\
P_5 = E_{[1,7]}+E_{[3,5]},\qquad & \qquad Q_5 = E_{[4,2]} + E_{[3,5]}, \\
P_6 = E_{[6,1]}+E_{[1,3]},\qquad & \qquad Q_6 = E_{[5,2]} + E_{[1,3]}, \\
\end{array}
$$
where a Cartan subalgebra $\hhh$ is spanned by  $P_0$ and $Q_0.$
The elements of the corresponding group $\HH = S^1\times S^1$ are of the form
$$
 g_{t,s} = \exp(t P_0 + s Q_0)= (\exp tP_0)(\exp sQ_0),\enskip t, s \in \RR .
$$
It is easy to show that the action of the element $g_{t,s}\in H$ on a point $p = (x_1,  x_2,  x_3, y_0, y_1,  y_2, y_3)\in S^6\subset
\RR ^7$ is given by
\begin{eqnarray}
g_{t,s}  p & = &  (x_1,x_2 \cos t -x_3 \sin t,  x_3 \cos t + x_2 \sin t, \label{eq:dejstvo}\\
& & \phantom{(x_1,,} y_0 \cos s + y_1 \sin s ,y_1 \cos s + y_0 \sin s, \nonumber \\
& & \phantom{(x_1,,} y_2 \cos (s-t)-y_3 \sin (s-t), y_3 \cos (s-t)+y_2 \sin (s-t)). \nonumber
\end{eqnarray}
Note that the action of  $H$ preserves  $x_1$ coordinate, so the orbit  ${\mathcal O}_p$ of any point $p$ belongs to the hyperplane $x_1 = const$ and therefore to the totally geodesic sphere $S^5\subset S^6.$ One can easily check that the orbit  ${\mathcal O}_p$ doesn't belong to any other hyperplane, i.e. the orbit ${\mathcal O}_p$ is not contained in a sphere of smaller dimension.

\noindent
The tangent space of the orbit ${\mathcal O}_p$ in the point $p = (x_1,  x_2,  x_3, y_0, y_1,  y_2, y_3)\in S^6$ is spanned by the vectors
\begin{eqnarray}
\overline{X}& = & \frac{d}{dt}(g_{t,s}p)|_{(t,s)=(0,0)}=
(0,-x_3,x_2,0,0,y_3,-y_2),\nonumber\\
\overline{Y}& = & \frac{d}{ds}(g_{t,s}p)|_{(t,s)=(0,0)}=
(0,0,0,-y_1,y_0,-y_3,y_2).\label{eq:baza}
\end{eqnarray}
These two vectors are linearly independent, i.e. the orbit  ${\mathcal O}_p$ is two-dimensional if the following conditions are satisfied
\begin{eqnarray}
\alpha & = & x_2^2 + x_3^2 + y_0^2 + y_1^2  \neq  0, \nonumber\\
\beta & = & x_2^2 + x_3^2 + y_2^2 + y_3^2  \neq  0, \label{eq:sfere}\\
\gamma & = & y_0^2 + y_1^2 + y_2^2 + y_3^2  \neq  0.\nonumber
\end{eqnarray}
Note that the corresponding equations represent three two-dimensional spheres on $S^6.$ In the sequel we consider the orbit of a point   $p\in S^6$ satisfying  relations $(\ref{eq:sfere}).$

\noindent
An orthonormal basis of tangent space  $T_p{\mathcal O}_p$ reads
\begin{eqnarray*}
X& =& \frac{1}{\sqrt{\beta}} (0,-x_3,x_2,0,0,y_3,-y_2),\\
Y& =& \frac{(0, -x_3(y_2^2 + y_3^2),x_2(y_2^2 + y_3^2), -\beta y_1, \beta y_0, y_3(x_2^2 + x_3^2), y_2(x_2^2 + x_3^2) )}{\sqrt{\beta}\sqrt{\alpha \beta - (x_2^2 + x_3^2)^2}}.
\end{eqnarray*}

\noindent
Using Lemma \ref{the:nezavodvekt} we get the slant angle of the orbit ${\mathcal O}_p$ in the point $p$ which, as we know, is constant along the orbit
\begin{equation}
\label{eq:theta}
cos\theta_p=|\langle
X,pY\rangle|=\frac{|x_3y_1y_2 - x_2y_0y_2 - x_2y_1y_3 -x_3y_0y_3|}{\sqrt{\alpha \beta - (x_2^2 + x_3^2)^2}}.
\end{equation}

\noindent
Now, we would like to know if we can find an orbit with slant angle $\theta$ for any given angle $\theta \in [0, \frac{\pi}{2}].$
From the formulas of the action (\ref{eq:dejstvo}) it follows that it is sufficient to  consider orbits of points satisfying $x_2 = 0 = y_0,$ i.e. $p =(x_1,  0,  x_3, 0, y_1,  y_2, y_3).$ Having in mind the slant angle (\ref{eq:theta}) we parameterize  such points $p$ in the following way
\begin{eqnarray*}
x3&=& \sqrt{1-x_1^2} \sin \theta \cos \varphi, \\
y1&=& \sqrt{1-x_1^2} \sin \varphi,\\
y2&=& \sqrt{1-x_1^2} \cos \theta \sin \phi \cos \varphi,\\
y3&=& \sqrt{1-x_1^2} \cos \theta \cos \phi \cos \varphi,
\end{eqnarray*}
where $x_1\in [-1,1], \enskip \theta, \varphi \in [-\frac{\pi}{2}, \frac{\pi}{2}], \enskip \phi \in [0, 2 \pi).$

\noindent
In the new coordinates the slant angle  (\ref{eq:theta}) becomes
\begin{equation}
\label{eq:th}
\cos \theta _p =
\frac{ \sin 2\varphi \, \sin 2\theta}{2\sqrt{4 \sin ^2 \varphi + \cos^2 \varphi \, \sin ^2 2\theta }} \sqrt{1-x_1^2} \, \sin \phi .
\end{equation}
Careful analysis of the above expression shows that the orbits ${\mathcal O}_p$ with minimal slant angle
$\theta = \arccos \frac{1}{3}$ correspond tho points
$$
p = \pm \frac{1}{\sqrt{3}}(0,0,1,0, 1,1,0).
$$
It is clear that there also exist points $p$ with totaly real orbits ${\mathcal O}_p,$ i.e. with  $\cos \theta = 0.$ Note  that there are no almost complex orbits.
By the previous consideration we prove the following theorem.
\begin{Theorem}
\label{th:orbits}
Let  $p\in S^6$ be a point satisfying relations  (\ref{eq:sfere}).
The orbit ${\mathcal O}_p$ under the action of Cartan subgroup $H\subset G_2$ on the sphere  $S^6$ is slant torus fully contained in totally geodesic $S^5\subset S^6.$ Its slant angle is given by the formula (\ref{eq:th}) and takes all values from the interval $[\arccos \frac{1}{3},\frac{\pi}{2}].$
\end{Theorem}

\noindent
Now, we analyze the geometry of the orbit ${\mathcal O}_p$ as a submanifold of sphere  $S^6$. Starting from the basis (\ref{eq:baza}) of $T_p{\mathcal O}_p$ one can calculate the induced connection and second fundamental form in the sphere  $S^6.$
Then, one can check that the Gausian curvature of the orbit  ${\mathcal O}_p$ vanishes. This also trivially follows from Gauss-Bonnet formula and the fact that Gausian curvature is constant along the orbit that is topologically a tori.

\noindent
One can check that the mean curvature vector of the orbit ${\mathcal O}_p$ in the point $p = (x_1,  0,  x_3, 0, y_1,  y_2, y_3)$ is given by the formula
$$
H(p) =
(2 x_1,0,\frac{N(x_3, y_1)}{D},0,\frac{N(y_1, x_3)}{D},y_2 (2-\frac{x_3^2+y_1^2}{D}),y_3(2-\frac{x3^2+y1^2}{D})),
$$
where
\begin{eqnarray*}
D &=& (y_1^2+y_2^2+y_3^2) x_3^2+y_1^2 (y_2^2+y_3^2),\\
N(y_1, x_3) &=& y_1 ((2 y_1^2+2 y_2^2+2 y_3^2-1) x_3^2+(2 y_1^2-1) (y_2^2+y_3^2)).
\end{eqnarray*}

\noindent
One can easily check that  $H\equiv 0$  for all points of the orbit  ${\mathcal O}_p,$ corresponding to the point
\begin{equation}
\label{eq:min}
  \quad p (\phi) = \frac{1}{\sqrt{3}}(0, 0, 1, 0,1, \cos \phi, \sin \phi), \enskip \phi \in [0,  2 \pi),
\end{equation}
i.e. exactly those orbits are minimal surfaces of $S^6.$

\noindent
From the formula (\ref{eq:th}) we get  the slant angle for all minimal orbits
$$
\theta  =
 \arccos \frac{\cos \phi}{3}.
$$
Therefore, we have the following theorem.
\begin{Theorem}
\label{th:minimal}
Let point $p\in S^6$ satisfies conditions (\ref{eq:sfere}). Its orbit  ${\mathcal O}_p$  under the action of Cartan subgroup $H\subset G_2$ is flat tori of $S^6.$ There exist one-parameter family of minimal orbits  ${\mathcal O}_p$ corresponding to the points (\ref{eq:min}).  For each  $ \theta \in [\arccos \frac{1}{3}, \frac{\pi}{2}]$ there exist exactly two minimal orbits  with slant angle $\theta .$
\end{Theorem}

\begin{Remark}
Totally real minimal orbits are obtained for $\phi = \frac{\pi}{2}, \frac{3\pi}{2}.$ They are  found in the classification of totally real minimal surfaces of constant curvature (see \cite{BVW97}, Theorem 6.5 (c)).
Minimal orbits with arbitrary slant angle are those from the classification of minimal, flat, slant surfaces of $S^6$ (see \cite{Xing}, Example 3.2).
\end{Remark}

\end{document}